\documentclass[12pt,a4paper]{amsart}
\usepackage[dvips]{graphicx}
\usepackage{amscd}
\usepackage{amsmath}
\usepackage{amssymb}
\usepackage{latexsym}
\usepackage[T1]{fontenc}
\usepackage[latin1]{inputenc}

\title[Brill--Noether theory on a graph]{
Brill-Noether theory of squarefree modules supported on a graph}

\author{Gunnar Fl{\o}ystad}
\address{Matematisk Institutt\\
         Johs. Brunsgt. 12\\
         5008 Bergen}
\email{gunnar@mi.uib.no}

\author{Henning Lohne}
\address{Matematisk Institutt\\
         Johs. Brunsgt. 12\\
         5008 Bergen}
\email{henning.lohne@math.uib.no}
\date{\today}

\begin{document}


\theoremstyle{plain}
\newtheorem{theorem}{Theorem}[section]
\newtheorem{corollary}[theorem]{Corollary}
\newtheorem*{main}{Main Theorem}
\newtheorem{lemma}[theorem]{Lemma}
\newtheorem{proposition}[theorem]{Proposition}
\newtheorem{conjecture}[theorem]{Conjecture}
\newtheorem*{pieri}{Pieri's rule}
\newtheorem*{theoremA}{Theorem A}
\newtheorem*{theoremB}{Theorem B}
\newtheorem{problem}{Problem}
\newtheorem*{Def2}{Definition}
\newtheorem*{Prop2}{Proposition}
\newtheorem*{Theorem2}{Theorem}

\theoremstyle{definition}
\newtheorem{definition}[theorem]{Definition}

\theoremstyle{remark}
\newtheorem{notation}[theorem]{Notation}
\newtheorem{remark}[theorem]{Remark}
\newtheorem{example}[theorem]{Example}
\newtheorem{claim}{Claim}


\newcommand{\psp}[1]{{{\bf P}^{#1}}}
\newcommand{\psr}[1]{{\bf P}(#1)}
\newcommand{\opw}{\op_{\psr{W}}}
\newcommand{\go}{\op}

\newcommand{\ini}[1]{\text{in}(#1)}
\newcommand{\gin}[1]{\text{gin}(#1)}
\newcommand{\kr}{{\Bbbk}}
\newcommand{\kk}{{\Bbbk}}
\newcommand{\pd}{\partial}
\newcommand{\vardel}{\partial}
\renewcommand{\tt}{{\bf t}}


\newcommand{\coh}{{{\text{{\rm coh}}}}}


\newcommand{\modv}[1]{{#1}\text{-{mod}}}
\newcommand{\modstab}[1]{{#1}-\underline{\text{mod}}}

\newcommand{\sut}{{}^{\tau}}
\newcommand{\sumit}{{}^{-\tau}}
\newcommand{\til}{\thicksim}

\newcommand{\totp}{\text{Tot}^{\prod}}
\newcommand{\dsum}{\bigoplus}
\newcommand{\dprod}{\prod}
\newcommand{\lsum}{\oplus}
\newcommand{\lprod}{\Pi}

\newcommand{\La}{{\Lambda}}
\newcommand{\lam}{{\lambda}}
\newcommand{\GL}{{GL}}

\newcommand{\sirstj}{\circledast}

\newcommand{\she}{\EuScript{S}\text{h}}
\newcommand{\cm}{\EuScript{CM}}
\newcommand{\cmd}{\EuScript{CM}^\dagger}
\newcommand{\cmri}{\EuScript{CM}^\circ}
\newcommand{\cler}{\EuScript{CL}}
\newcommand{\clerd}{\EuScript{CL}^\dagger}
\newcommand{\clerri}{\EuScript{CL}^\circ}
\newcommand{\gor}{\EuScript{G}}
\newcommand{\gF}{\mathcal{F}}
\newcommand{\gG}{\mathcal{G}}
\newcommand{\gM}{\mathcal{M}}
\newcommand{\gE}{\mathcal{E}}
\newcommand{\gD}{\mathcal{D}}
\newcommand{\gI}{\mathcal{I}}
\newcommand{\gP}{\mathcal{P}}
\newcommand{\gK}{\mathcal{K}}
\newcommand{\gL}{\mathcal{L}}
\newcommand{\gS}{\mathcal{S}}
\newcommand{\gC}{\mathcal{C}}
\newcommand{\gO}{\mathcal{O}}
\newcommand{\gJ}{\mathcal{J}}
\newcommand{\gU}{\mathcal{U}}
\newcommand{\mm}{\mathfrak{m}}

\newcommand{\dlim} {\varinjlim}
\newcommand{\ilim} {\varprojlim}

\newcommand{\CM}{\text{CM}}
\newcommand{\Mon}{\text{Mon}}


\newcommand{\Kom}{\text{Kom}}


\newcommand{\EH}{{\mathbf H}}
\newcommand{\res}{\text{res}}
\newcommand{\Hom}{\text{Hom}}
\newcommand{\inhom}{{\underline{\text{Hom}}}}
\newcommand{\Ext}{\text{Ext}}
\newcommand{\Tor}{\text{Tor}}
\newcommand{\ghom}{\mathcal{H}om}
\newcommand{\gext}{\mathcal{E}xt}
\newcommand{\id}{\text{{id}}}
\newcommand{\im}{\text{im}\,}
\newcommand{\codim} {\text{codim}\,}
\newcommand{\resol}{\text{resol}\,}
\newcommand{\rank}{\text{rank}\,}
\newcommand{\lpd}{\text{lpd}\,}
\newcommand{\coker}{\text{coker}\,}
\newcommand{\supp}{\text{supp}\,}
\newcommand{\Ad}{A_\cdot}
\newcommand{\Bd}{B_\cdot}
\newcommand{\Gd}{G_\cdot}


\newcommand{\sus}{\subseteq}
\newcommand{\sups}{\supseteq}
\newcommand{\pil}{\rightarrow}
\newcommand{\vpil}{\leftarrow}
\newcommand{\lvpil}{\longleftarrow}
\newcommand{\rpil}{\leftarrow}
\newcommand{\lpil}{\longrightarrow}
\newcommand{\inpil}{\hookrightarrow}
\newcommand{\pils}{\twoheadrightarrow}
\newcommand{\projpil}{\dashrightarrow}
\newcommand{\dotpil}{\dashrightarrow}
\newcommand{\adj}[2]{\overset{#1}{\underset{#2}{\rightleftarrows}}}
\newcommand{\mto}[1]{\stackrel{#1}\longrightarrow}
\newcommand{\vmto}[1]{\overset{\tiny{#1}}{\longleftarrow}}
\newcommand{\mtoelm}[1]{\stackrel{#1}\mapsto}

\newcommand{\eqv}{\Leftrightarrow}
\newcommand{\impl}{\Rightarrow}

\newcommand{\iso}{\cong}
\newcommand{\te}{\otimes}
\newcommand{\tek}{\te_\kr}
\newcommand{\sqte}{\te}
\newcommand{\into}[1]{\hookrightarrow{#1}}
\newcommand{\ekv}{\Leftrightarrow}
\newcommand{\equi}{\simeq}
\newcommand{\isopil}{\overset{\cong}{\lpil}}
\newcommand{\equipil}{\overset{\equi}{\lpil}}
\newcommand{\ispil}{\isopil}
\newcommand{\vvi}{\langle}
\newcommand{\hvi}{\rangle}
\newcommand{\susneq}{\subsetneq}
\newcommand{\sgn}{\text{sign}}
\newcommand{\Spec}{Spec\,}


\newcommand{\xd}{\check{x}}
\newcommand{\ortog}{\bot}
\newcommand{\tL}{\tilde{L}}
\newcommand{\tM}{\tilde{M}}
\newcommand{\tH}{\tilde{H}}
\newcommand{\tvH}{\widetilde{H}}
\newcommand{\tvh}{\widetilde{h}}
\newcommand{\tV}{\tilde{V}}
\newcommand{\tS}{\tilde{S}}
\newcommand{\tT}{\tilde{T}}
\newcommand{\tR}{\tilde{R}}
\newcommand{\tf}{\tilde{f}}
\newcommand{\ts}{\tilde{s}}
\newcommand{\tp}{\tilde{p}}
\newcommand{\tr}{\tilde{r}}
\newcommand{\tfst}{\tilde{f}_*}
\newcommand{\empt}{\emptyset}
\newcommand{\bfa}{{\bf a}}
\newcommand{\bfb}{{\bf b}}
\newcommand{\bfd}{{\bf d}}
\newcommand{\bfe}{{\bf e}}
\newcommand{\bfp}{{\bf p}}
\newcommand{\bfc}{{\bf c}}
\newcommand{\bfl}{{\bf \ell}}
\newcommand{\bfz}{{\bf z}}
\newcommand{\bfj}{{\bf j}}
\newcommand{\ubfd}{\underline{\bfd}}
\newcommand{\la}{\lambda}
\newcommand{\bfen}{{\mathbf 1}}
\newcommand{\ep}{\epsilon}
\newcommand{\en}{r}
\newcommand{\tu}{s}
\newcommand{\integ}{{int}}
\newcommand{\module}{{mod}}
\newcommand{\inc}{{deg}}
\newcommand{\dec}{{root}}

\newcommand{\ome}{\omega_E}

\newcommand{\bevis}{{\bf Proof. }}
\newcommand{\demofin}{\qed \vskip 3.5mm}
\newcommand{\nyp}[1]{\noindent {\bf (#1)}}
\newcommand{\demo}{{\it Proof. }}
\newcommand{\demodone}{\demofin}
\newcommand{\parg}{{\vskip 2mm \addtocounter{theorem}{1}  
                   \noindent {\bf \thetheorem .} \hskip 1.5mm }}

\newcommand{\lcm}{{\text{lcm}}}


\newcommand{\dl}{\Delta}
\newcommand{\cdel}{{C\Delta}}
\newcommand{\cdelp}{{C\Delta^{\prime}}}
\newcommand{\dlst}{\Delta^*}
\newcommand{\Sdl}{{\mathcal S}_{\dl}}
\newcommand{\lk}{\text{lk}}
\newcommand{\lkd}{\lk_\Delta}
\newcommand{\lkp}[2]{\lk_{#1} {#2}}
\newcommand{\del}{\Delta}
\newcommand{\delr}{\Delta_{-R}}
\newcommand{\dd}{{\dim \del}}

\renewcommand{\aa}{{\bf a}}
\newcommand{\bb}{{\bf b}}
\newcommand{\cc}{{\bf c}}
\newcommand{\xx}{{\bf x}}
\newcommand{\yy}{{\bf y}}
\newcommand{\zz}{{\bf z}}
\newcommand{\mv}{{\xx^{\aa_v}}}
\newcommand{\mF}{{\xx^{\aa_F}}}

\newcommand{\proj}[1]{{\mathbb P}^{#1}}
\newcommand{\hele}{{\mathbb Z}}
\newcommand{\nat}{{\mathbb N}}
\newcommand{\rat}{{\mathbb Q}}

\newcommand{\pnm}{{\bf P}^{n-1}}
\newcommand{\opnm}{{\go_{\pnm}}}
\newcommand{\op}[1]{\gO_{\proj{#1}}}
\newcommand{\ompn}{\Omega_{\proj{n}}}
\newcommand{\opn}{\op{n}}
\newcommand{\opm}{\op{m}}
\newcommand{\ompnm}{\omega_{\pnm}}

\newcommand{\dt}{{\displaystyle \cdot}}
\newcommand{\st}{\hskip 0.5mm {}^{\rule{0.4pt}{1.5mm}}}              
\newcommand{\disk}{\scriptscriptstyle{\bullet}}

\newcommand{\cF}{F_\dt}
\newcommand{\Fd}{F_{\disk}}
\newcommand{\pol}{f}

\newcommand{\disc}{\circle*{5}}

\newcommand{\Dab}{{\mathbb D}(\bfa, \bfb)}
\newcommand{\Ddab}{B(\bfa, \bfb)}
\newcommand{\Tab}{{\mathbb T}(\bfa, \bfb)}
\newcommand{\Cab}{C(\bfa, \bfb)}
\newcommand{\Pab}{B(\bfa, \bfb)}
\newcommand{\Lhkab}{L^{HK}(\bfa, \bfb)}
\newcommand{\Lab}{L(\bfa, \bfb)}
\newcommand{\Sigab}{\Sigma(\bfa, \bfb)}
\newcommand{\hlow}{h_{low}}
\newcommand{\hup}{h_{up}}
\newcommand{\up}[2]{\hup}
\newcommand{\facet}[2]{{\bf facet}(#1, #2)}
\newcommand{\BS}{Boij-S\"oderberg}
\newcommand{\bstar}{{\mathbb B}^*}
\newcommand{\Symm}{\mbox{Symm}}
\newcommand{\zinc}[1]{{\mathbb Z}^{#1}_\inc}
\newcommand{\zdec}[1]{{\mathbb Z}^{#1}_\dec}

\def\CC{{\mathbb C}}
\def\GG{{\mathbb G}}
\def\ZZ{{\mathbb Z}}
\def\NN{{\mathbb N}}
\def\RR{{\mathbb R}}
\def\OO{{\mathbb O}}
\def\QQ{{\mathbb Q}}
\def\VV{{\mathbb V}}
\def\PP{{\mathbb P}}
\def\EE{{\mathbb E}}
\def\FF{{\mathbb F}}
\def\AA{{\mathbb A}}

\begin{abstract}
We investigate the analogy between squarefree Cohen-Macaulay modules
supported on a graph and line bundles on a curve. We prove a Riemann--Roch 
theorem, we study the Jacobian and gonality of a graph, and
we prove Clifford's theorem.
\end{abstract}

\bibliographystyle{plain}

\maketitle

\section{Introduction}
Let $S=k[x_1,...,x_n]$ be the polynomial ring. K.\ Yanagawa (\cite{Ya1} and 
\cite{Ya2}) introduced the notion of squarefree $\mathbf{Z}^n$-graded 
$S$-modules.
Squarefree modules is a generalization of squarefree monomial ideals. These
modules are still of a combinatorial nature, as they have support on a 
simplicial complex.
Such modules can be thought of as combinatorial analogues to coherent sheaves,
or vector bundles if they are Cohen--Macaulay.

In this paper we study CM squarefree 
modules whose support is a graph.
These may be thought of as analogues of line bundles on smooth projective
curves. Inspired by this we show how 
results from the theory of line bundles on 
curves have analogues in the theory of squarefree CM modules on graphs.
We define their degree and global 
sections, and we show an analog of the Riemann--Roch theorem 
(Theorem \ref{RR-Thm}) for simplical graphs. There is another combinatorial 
Riemann--Roch theorem for graphs, shown by M.\ Baker and S.\ Norine in 
\cite{BN}. Our setting is however different, and our result should be thought
of as an algebraic combinatorial analogue. Furthermore, we study squarefree 
modules of multi-degree $(0,\dots,0)$ with support on a graph. 
This family is an analogue to the Jacobian variety of a curve. We show that 
the dimension of the family of isomorphism classes of such modules equals the
genus of the simplicial graph, similar as for the moduli space of line bundles
on a curve. 

Then we consider results from Brill--Noether theory, and must then
limit ourselves to two-connected graphs. 
We define effective and special squarefree modules on such 
graphs, and we define an analogue of gonality. 
We show that the gonality has the same upper bound as in the classical case.
In the end we 
prove Clifford's theorem for graphs. In Remark \ref{Graph_curves}, 
we briefly discuss the possible connections between these numerical invariants 
and the resolution of graph curves. A similar study has been done by D.\ Bayer
and D.\ Eisenbud in \cite{BE}.

\medskip
The paper is organized as follows: Section \ref{prelim} contains 
preliminaries about squarefree $S$-modules. We also recall a combinatorial 
description of the canonical module associated to the Stanley-Reisner
ring of a simplicial complex.

In Section \ref{RR}, we describe the analogy between squarefree modules and 
vector bundles. The main theorem of this section is the analogue of the 
Riemann--Roch theorem. We also give upper and lower bounds for the degree of 
an indecomposable squarefree module. We also  
show the existence of indecomposable squarefree modules of degree $d$ 
with a one-dimensional space of sections, when $0 \leq d \leq g$, the 
genus of a graph, analogous the case of line bundles on a curve.

In Section \ref{Jacobian}, we study squarefree modules of multi-degree 
$(0,\dots,0)$ with support on a simplicial graph. We show that the dimension 
of the family of isomorphism classes of such modules are the genus of the 
graph, similar as for the moduli space of line bundles on a curve.

In Section \ref{Go}, we define 
effective and special squarefree modules on two-connected graphs, and study the
analogue of gonality. The main result is that the gonality has the
same upper bound as in the classical case.
In the last Section \ref{Cl} we prove Clifford's theorem

\medskip
\noindent
\textbf{Acknowledgments.} The second author thanks Professor 
Frank-Olaf Schreyer for helpful comments during his stay at the University of 
Saarland. 
We are also grateful to the referee for several suggestions to improve
the paper, and for 
pointing out gaps in some of the original proofs.

\section{Preliminaries}\label{prelim}

The most natural squarefree modules are squarefree monomial ideals and the corresponding quotient rings. These objects are given by the combinatorial structure of simplicial complexes and Stanley--Reisener rings. We recall the following definitions.

Let $S:=k[x_1,x_2,\dots,x_n]$ be the polynomial ring over some field $k$ and let $[n]=\{1,2,\dots,n\}$. A subset $F\subseteq [n]$ is called a face. A simplicial complex is a collection of faces $\Delta$, such that if $F\in \Delta$ and $G\subseteq F$, then $G\in \Delta$. The Stanley--Reisner ideal of the simplicial complex $\Delta$ is the squarefree monomial ideal $I_\Delta=<x^\sigma\,|\, \sigma\not\in \Delta >$ generated by monomials corresponding to non-faces $\sigma$ of $\Delta$. The Stanley--Reisner ring of $\Delta$ is the quotient ring $k[\Delta]:=S/{I_{\Delta}}$

\begin{definition}[cf. \cite{Ya2} Definition 2.1]
A finitely generated $\mathbf{Z}^n$-graded $S$-module $M=\bigoplus_{\mathbf{a}\in\mathbf{Z}^n} M_\mathbf{a}$ is said to be squarefree, if the following hold
\begin{itemize}
\item[i)] $M$ is $\mathbf{N}^n$-graded, that is $M_\mathbf{a}=0$ for $\mathbf{a}\not\in \mathbf{N}^n$.
\item[ii)] The multiplication map $M_\mathbf{a}\overset{\cdot x_i}{\rightarrow} M_{\mathbf{a}+\mathbf{e}_i}$ is bijective for all $\mathbf{a}\in\mathbf{N}^n$ and all
 $i$ in  Supp$(\mathbf{a})$.
\end{itemize}
\end{definition}

\begin{theorem}[cf. \cite{Ya1} Theorem 2.6]
If $M$ is a squarefree $S$-module, then so is $\emph{\text{Ext}}_S^i(M,\omega_S)$ for all $i$.
\end{theorem}

Let $M$ be a $\mathbf{Z}^n$-graded $S$-module. The Matlis dual of $M$ is the $\mathbf{Z}^n$-graded $S$-module $M^\vee=\text{Hom}_k(M,k)$. This means that $$(M^\vee)_{-\mathbf{a}}=\text{Hom}_k(M_\mathbf{a},k),$$ and the multiplication map $(M^\vee)_{-\mathbf{a}}\overset{\cdot \mathbf{x}^\mathbf{b}}{\rightarrow}(M^\vee)_{-\mathbf{a}+\mathbf{b}}$ is the transpose of the multiplication map $M_{\mathbf{a}-\mathbf{b}}\overset{\cdot \mathbf{x}^\mathbf{b}}{\rightarrow}M_\mathbf{a}$. See \cite{Mi} and \cite[Def. 11.15]{MS}.

Let $H_m^i(M)$ denote the local cohomology modules of $M$. Then the following holds.

\begin{theorem}[Local duality] \label{PreTheLokdul}
For all finitely generated $\mathbf{Z}^n$-graded $S$-modules $M$ and all integers $i$ there exist natural homogeneous isomorphisms
$$(H_m^i(M))^\vee \cong \emph{Ext}_S^{n-i}(M,\omega_S).$$
\end{theorem}

\begin{proof}
See \cite[Thm. 3.6.19]{BH} or \cite[Cor. 6.1]{Mi}.
\end{proof}

The local cohomology modules can be calculated by the cocomplex
\begin{equation} \label{PreLigKcx}
K^\bullet:0\rightarrow K^0\rightarrow K^1\rightarrow \cdots \rightarrow 
K^n\rightarrow 0,
\end{equation}
where 
$$K^i=\underset{|\sigma|=i}{\underset{\sigma\subseteq [n]}{\bigoplus}} 
M_{\mathbf{x}^\sigma}$$
is the direct sum of the module localized in the squarefree monomials of degree
 $i$, and the differential is given on each component as
$$\text{sign}(j,\sigma\cup j)\cdot \text{nat}: M_{\mathbf{x}^\sigma}\rightarrow 
M_{\mathbf{x}^{\sigma\cup j}}$$
if $j\not\in \sigma$, where $\text{sign}(j,\sigma)=(-1)^{\#\{i\in\sigma|i<j\}}$. 
If $M$ is squarefree, then the cohomology $H^i(K^\bullet)\cong H_m^i(M)$ is easy 
to calculate. As in \cite{Ya1} we have that if $M$ is squarefree, then 
$$(M_{\mathbf{x}^\sigma})_{-\tau}=\left\{
\begin{array}{ll}
M_\sigma & \text{if }\tau\subseteq \sigma \\
0 & \text{otherwise},
\end{array}
\right.
$$
and the natural map $(M_{\mathbf{x}^\sigma})_{-\tau}\rightarrow 
(M_{\mathbf{x}^{\sigma\cup \{i\}}})_{-\tau}$ corresponds to the map 
$M_\sigma\overset{\cdot x_i}{\rightarrow} M_{\sigma\cup\{i\}}$.

\begin{definition}
Let $M$ be a squarefree Cohen--Macaulay $S$-module of dimension $d$. 
The canonical module of $M$ is defined as the squarefree $S$-module
$$\omega_M=\mathrm{Ext}_S^{n-d}(M,\omega_S).$$
\end{definition}

\begin{proposition}\label{omegadim} If $M$ is a squarefree Cohen--Macaulay 
module of dimension $d$ and if $\tau\subseteq [n]$, then
\[ \emph{dim}_k (\omega_M)_\tau =\sum_{\tau\subseteq \sigma\subseteq [n]} 
(-1)^{d-|\sigma|}\emph{dim}_k M_\sigma . \]
\end{proposition}

\begin{proof}
If $M$ is Cohen--Macaulay of dimension $d$, then the local cohomology groups 
$H_m^i(M)=0$ for $i\neq d$. So the cocomplex $K^\bullet$ restricted to degree 
$-\tau$ has only cohomology $(H_m^d(M))_{-\tau}$, which is isomorphic to 
$\text{Hom}_k((\omega_M)_\tau,k)$ by local duality. The result now follows since 
$$(K^i)_{-\tau}=\underset{\tau\subseteq \sigma}{\underset{|\sigma|=i}{\bigoplus}}
M_\sigma .$$
\end{proof}

\subsection{The canonical module of $k[\Delta]$}
Assume that $\Delta$ is a Cohen--Macaulay $(d-1)$-dimensional simplicial 
complex, i.e., the Stanley--Reisner ring $k[\Delta]$ is Cohen--Macaulay of 
dimension $d$. The Matlis dual of the cocomplex $K^\bullet$ above restricted to 
positive degrees gives a nice $\mathbf{Z}^n$-graded description of 
$\omega_{k[\Delta]}$. The details are written out in \cite[Theorem 5.7.3]{BH}, 
here we only give the result.

For each $i=0,\dots,d$, let $G_i$ be the direct sum of the $k[\Delta]$-modules 
$k[X_1,\dots,X_n]/I_F$ where $F\in\Delta$, $|F|=i$ and $I_F=(X_j\,|\, j 
\not\in F)$. Let $\varphi_i:G_i\rightarrow G_{i-1}$ be the map which is 
$$(-1)^{j-1}\text{nat}:k[X_1,\dots,X_n]/I_F\rightarrow k[X_1,\dots,X_n]/I_{F'}$$
if $F=\{v_{i_1},\dots,v_{i_r}\}$ and $F'=\{v_{i_1},\dots,\widehat{v_{i_j}},
\dots,v_{i_r}\}$, and zero otherwise. Then the complex
\begin{equation} \label{PreLigGcx}
G_{\bullet}:0\rightarrow G_d\rightarrow G_{d-1}\rightarrow\cdots\rightarrow 
G_1\rightarrow G_0\rightarrow 0 
\end{equation}
is the Matlis dual of $K^\bullet$ restricted to positive degrees.
When $\Delta$ is Cohen--Macaulay we obtain an exact sequence of 
$\mathbf{Z}^n$-graded $k[\Delta]$-modules
$$0\rightarrow \omega_{k[\Delta]}\rightarrow G_d\rightarrow G_{d-1}\rightarrow 
\cdots \rightarrow G_1\rightarrow G_0\rightarrow 0.$$

The description of the canonical module in the long exact sequence above gives rise to a nice description of its squarefree grades. Let $\text{st}_\Delta F$ denote the set of faces of $\Delta$ containing $F$, and let $\Delta-F$ denote $\Delta\setminus \text{st}_\Delta F$. We have that
$$(\omega_{k[\Delta]})_0\cong \widetilde{H}_{d-1}(\Delta;k) \text{ and}$$
$$(\omega_{k[\Delta]})_F\cong \widetilde{H}_{d-1}(\Delta,\Delta - F;k)\text{ for faces } \emptyset\neq F\in \Delta, $$
and the multiplication map
$$\widetilde{H}_{d-1}(\Delta,\Delta - F;k)\overset{\cdot x_i}{\rightarrow} \widetilde{H}_{d-1}(\Delta,\Delta - (F\cup\{i\});k)$$
is the natural map.
\\

\begin{remark}
For any face $F$, the homology groups of the chain complex above, restricted to degree $F$, can also be described in the more common way using the link of a face. Recall that the link of $F$ in $\Delta$ is defined as $\mathrm{lk}_F \Delta := \left\{G\in \Delta \, \middle| \, F\cap G = \emptyset \text{ and } F\cup G \in \Delta \right\}$. The homology groups above can now be described as $\widetilde{H}_i(\Delta,\Delta-F;k)\cong \widetilde{H}_{i-|F|}(\text{lk}_\Delta F;k).$ This is because the chain complexes are the same. The two different descriptions both have their advantages. The first gives a natural description of the multiplication map, while the second is a combinatorial description.
\\
\end{remark}

We use these homology groups to give a characterisation of $2$-CM simplicial complexes. Recall that a Cohen--Macaulay simplicial complex $\Delta$ is said to be doubly Cohen--Macaulay or $2$-CM if $\Delta -\{p\}$ is Cohen--Macaulay of the same dimension as $\Delta$ for all vertices $p\in\Delta$. The following results might be well-known to specialists.
However we give proofs for the reader's convenience.

\begin{proposition}
Let $\Delta$ be a simplicial complex of dimension $d-1$. Then the following are equivalent.
\begin{itemize}
\item[1.] $\Delta$ is $2$-CM.
\item[2.] $\widetilde{H}_i(\Delta;k)=0$ for $0\le i\le d-2$ and $\widetilde{H}_i(\Delta-F;k)=0$ for $0\le i\le d-2$ and any face $F\in \Delta$.
\end{itemize}
\end{proposition}

\begin{proof}
The homology groups $\widetilde{H}_i(\Delta;k)$ and 
$\widetilde{H}_i(\Delta,\Delta - F;k)$ are the generators for the module 
$\mathrm{Ext}_S^{n-i-1}(k[\Delta], \omega_S)$ by the complex (\ref{PreLigGcx}), 
and by local duality Theorem \ref{PreTheLokdul}. 
Because of the long exact sequence
\begin{eqnarray*}
\cdots & \rightarrow &  \widetilde{H}_i(\Delta;k)\rightarrow
\widetilde{H}_i(\Delta,\Delta - F;k)\rightarrow 
\widetilde{H}_{i-1}(\Delta - F;k) \\
& \rightarrow & \widetilde{H}_{i-1}(\Delta;k)\rightarrow \cdots
\end{eqnarray*}
 we see that 
$\Delta$ is CM if and only if $\widetilde{H}_i(\Delta;k)=0$ for 
$0\le i\le d-2$ and $\widetilde{H}_{i}(\Delta - F;k)=0$ for $0\le i\le d-3$ and
 all faces $\emptyset\neq F\in \Delta$. So it is enough to show that $\Delta$ 
is $2$-CM if and only if $\Delta$ is CM and that 
$\widetilde{H}_{d-2}(\Delta-F;k)=0$ for any face $F\in \Delta$.

Suppose that  $\Delta$ is $2$-CM. We will show that $\widetilde{H}_{d-2}(\Delta - F;k)=0$ for all faces $F\in\Delta$ using induction on the dimension of $F$. If $F=\{p\}$ is just a vertex, then $\widetilde{H}_{d-2}(\Delta-\{p\};k)=0$ 
since $\Delta - p$ is CM. Now suppose that $G=F\cup \{p\}$ is a face. 
Let $A=\Delta-F$ and $B=\Delta - \{p\}$. Then $A\cup B = \Delta - G$ and $A\cap B = (\Delta - \{p\}) - F$. We therefore have a Mayer--Vietoris sequence
\begin{eqnarray*}
\cdots & \rightarrow & \widetilde{H}_{d-2}(\Delta - F;k)\oplus
\widetilde{H}_{d-2}(\Delta - \{p\};k)\rightarrow \widetilde{H}_{d-2}(\Delta-G;k)
\\ & \rightarrow & \widetilde{H}_{d-3}((\Delta - \{p\}) -F;k)\rightarrow \cdots.
\end{eqnarray*}
The homology groups on the left side are trivial because of the induction 
hypothesis and the homology group on the right side is trivial because 
$\Delta - \{p\}$ is CM of dimension $d-1$ so it follows that 
$\widetilde{H}_{d-2}(\Delta - G;k)=0$.

Suppose that $\Delta$ is CM and that $\widetilde{H}_{d-2}(\Delta - F;k)=0$ 
for all faces $F\in\Delta$. We have that $\widetilde{H}_i(\Delta - \{p\};k)=0$ 
for $0\le i \le d-2$. Let $p\in \Delta$ be a vertex and $F\in \Delta$ a face 
not containing $p$. As above, we get a Mayer--Vietoris sequence
\begin{eqnarray*}
\cdots \rightarrow \widetilde{H}_{i+1}(\Delta-G;k) & \rightarrow & 
\widetilde{H}_{i}((\Delta-\{p\})-F;k) \\
& \rightarrow & 
\widetilde{H}_{i}(\Delta-\{p\};k)\oplus 
\widetilde{H}_{i}(\Delta-F;k)\rightarrow \cdots,
\end{eqnarray*} where $G$ is the face 
$F \cup \{p\}$ if this is a face in $\Delta$, otherwise, we just replace $\Delta-G$ above with $\Delta$. For $0\le i \le d-3$, the homology group on the left side is trivial because of the assumption and the fact that $\Delta$ is CM, and the homology groups on the right side are trivial because $\Delta$ is CM, so $\widetilde{H}_{i}((\Delta-\{p\})-F;k)=0$ for $0\le i\le d-3$ and it follows that $\Delta - p$ is CM for all vertices $p\in \Delta$, hence $\Delta$ is $2$-CM. 
\end{proof}

This characterisation of $2$-CM simplicial complexes has the following 
corollary, which we will use later. For more details see \cite{Ba}, \cite{BH}
and \cite{St}.

\begin{corollary}\label{2CM}
Let $\Delta$ be a CM simplicial complex of dimension $d-1$. Then $\Delta$ is 
$2$-CM if and only if $\omega_{k[\Delta]}$ is generated by 
$(\omega_{k[\Delta]})_0\cong \widetilde{H}_{d-1}(\Delta;k)$.
\end{corollary}

\begin{proof}
Consider the long exact sequence
$$0\rightarrow \widetilde{H}_{d-1}(\Delta - F)\rightarrow \widetilde{H}_{d-1}(\Delta)\overset{\cdot\mathbf{x}^F}{\rightarrow} \widetilde{H}_{d-1}(\Delta,\Delta - F)\rightarrow \widetilde{H}_{d-2}(\Delta - F)\rightarrow 0.$$
$\omega_{k[\Delta]}$ is generated by $\widetilde{H}_{d-1}(\Delta;k)$ if and only if $\cdot\mathbf{x}^F$ is surjective for every $F$. That is, if and only if $\widetilde{H}_{d-2}(\Delta - F)=0$ for every face $F$.

\end{proof}

\section{The Riemann--Roch theorem}\label{RR}

The annihilator of a squarefree $S$-module $M$ is a squarefree monomial ideal. Since squarefree monomial ideals corresponds to simplicial complexes, $M$ can be considered as a module over the Stanley--Reisner ring $k[\Delta]=k[x_1,\dots,x_n]/\text{ann}(M)$. The study of squarefree Cohen--Macaulay modules are therefore the same as the study of maximal Cohen--Macaulay modules over Stanley--Reisner rings $k[\Delta]$, with support on all of $k[\Delta]$. 
A squarefree $S$-module can be described as follows.
 For each face $A\subseteq [n]$ we have a $k$-vector space $M_A$, 
and for each face $A\subseteq B$ we have a $k$-linear map 
$\varphi_{AB}:M_A\rightarrow M_B$ satisfying $\varphi_{AA}=id_{M_A}$ and 
if $A\subseteq B\subseteq C$, then $\varphi_{BC}\circ \varphi_{AB}=\varphi_{AC}$.

In the rest of this paper we study squarefree Cohen--Macaulay modules
with support on a graph. We always assume the graph to be 
{\it connected, simple and
without isolated vertices}, so
it is a CM one-dimensional simplicial complex. 
First note the following simple criterion.

\begin{lemma} \label{RRLemCM} A squarefree module $M$ with support a graph is
CM if and only if (the $e$ denote edges):
\begin{itemize}
\item[1.] For each vertex $v$ the map 
\[M_v \mto{\oplus \phi_{ve}} \oplus _{v\in e} M_e \]
is injective.
\item[2.] The following sequence is left exact:
\[ 0 \pil M_{\emptyset} \mto{\oplus \phi_{\emptyset v}}
\underset{\text{vertices } v}{\oplus}  M_v \mto{\oplus \phi_{\emptyset v}}
\underset{\text{edges } e}{\oplus} M_e. \]
\end{itemize}
\end{lemma}

\begin{proof} That $M$ is CM means that the complex (\ref{PreLigKcx})
has cokernel only in cohomological degree $2$. When considering it
in degree $\{v \}$ we get 1. above, and when considering it in degree
$\emptyset$ we get 2.
\end{proof}

Such a module $M$ on a graph $G$ gives rise to a
sheaf $\widetilde{M}$ on $\mathrm{Proj}\,\,k[G]$ and there is a natural
graded isomorphism between $M$ and the graded module of global sections 
$\bigoplus_\tau \Gamma(\mathrm{Proj}\,\,k[G],\widetilde{M}(\tau))$ 
\cite[Appendix A4]{Ei}. We shall consider such modules with 
$M_e$ one-dimensional for all edges $e$ in the graph.
 In this case we say that $M$ is locally of rank $1$. 
Such a module may be seen as the analog of a line bundle on a curve.

Inspired by this,
 we define $l(M):=\text{dim}_k M_\emptyset$ as an analogue of the global 
sections. Furthermore, we define the multi-degree of $M$ as the vector 
$\mathbf{d}\in \mathbf{Z}^n$, such that $d_i=\text{dim}_k M_{\{i\}} -1$, 
and the degree of $M$ as $\text{deg}(M):=\sum_i d_i$. This definition of the
degree of a module is an analogy to the degree of a line bundle. 

Let $V$ be the vertices of the graph $G$ and $E$ its set of its edges.
We define the genus of $G$ as:
$$g(G ):=l(\omega_{k[G]})=\text{dim}_k (\omega_{k[G]})_0=\text{dim}_k \widetilde{H}_1(G ;k)=1-|V|+|E|.$$

\begin{theorem}[Riemann-Roch]\label{RR-Thm}
Let $G$ be a graph, and $M$ a squarefree Cohen--Macaulay module with support 
on $G$, locally of rank $1$. Then the following formula holds.
$$l(M)-l(\omega_M)=1+\emph{deg}(M)-g.$$
\end{theorem}
\begin{proof} Using Proposition \ref{omegadim}, we can calculate the left hand side as $$\sum_{|\sigma|=1} M_\sigma -\sum_{|\tau|=2} M_\tau=(\text{deg}(M)+|V|)-|E|
=\text{deg}(M)+1-g.$$
\end{proof}

Using the Riemann--Roch formula, it is easy to see that 
$\text{deg}(\omega_M)=2g-2-\text{deg}(M)$. Since the degree is obviously
bounded below by $\text{deg}(M)\ge -|V|$ and $\omega_M$ is a squarefree 
Cohen--Macaulay module whenever $M$ is, we get that the degree is also bounded
above by $\text{deg}(M)\le 2g-2+|V|$. For a given   graph, it is possible to
give a better bound for the degree.

If $\tau$ 
is a subset of the vertices, the {\it restriction $G|_\tau$} is the graph
whose vertices are $\tau$ and whose edges are the edges in $G$ with
both endpoints in $\tau$.

\begin{proposition}
Let $G$ be a graph, and $M$ an indecomposable CM squarefree module on $G$,
locally of rank $1$. Then
$$-s\le\emph{deg}(M)\le 2g-2+s,$$
where 
$$s=\emph{max}\left\{|\sigma| \, \middle| \, \sigma\subseteq V \text{, }
\emph{dim}(G|_\sigma)=0\text{ and } G|_{\sigma^c} \text{ is connected}\right\}.$$
Furthermore, for any $-s\le i \le 2g-2+s$ there exists such a module $M$ with $\mathrm{deg}(M)=i$.
\end{proposition}

\begin{proof}
Since the upper bound is the dual of the lower bound, it is enough to show that the lower bound holds. Suppose that $\text{deg}(M)=-k$. Then there are at least $k$ vertices in $G$ where $M_v=0$. Let $\tau\subseteq [n]$ be the subset that corresponds to these vertices. Since $M$ is indecomposable, $G|_{\tau^c}$ must be connected and there can not be any edges in $G|_\tau$, so $\text{dim}(G|_\tau)=0$. It therefore follows from the definition of $s$ in the proposition that $-s\le -k$.

For $0\le i \le s$, let $\sigma\subseteq [n]$ be a subset such that $|\sigma|=i$, $\text{dim}(G|_\sigma)=0$ and $G|_{\sigma^c}$ is connected. Then we can construct an indecomposable module $M$ with $\text{deg}(M)=-i$. Let $M$ be the module where $M_v=0$ for all vertices $v\in \sigma$, $M_v = k 1_v $ for all vertices 
$v\in \sigma^c$ and $M_e = k 1_e$ for all edges, with multiplication map 
$\varphi_{v e}(1_v) = 1_e$ for all vertices $v\in \sigma^c$ and
edges $e$ containing $v$. 

If $M$ is decomposable, isomorphic to $M^\prime \oplus M^{\prime \prime}$, let 
$V^\prime$ be the vertices in $\sigma^c$ at which $M^\prime_v$ is nonzero
and $V^{\prime \prime}$ the vertices of $\sigma^c$ at which $M^{\prime \prime}$
is nonzero. Then $V^\prime \cup V^{\prime \prime}$ is a partition of $\sigma^c$. 
The edges of $G$ may be partitioned into $E^\prime \cup E^{\prime \prime}$
where $E^\prime$ are the edges in the support of $M^\prime$ and correspondingly
for $E^{\prime \prime}$. Note that since the maps $\phi_{ve}$ are all nonzero
when $v \in \sigma^c$ and $v \in e$, an edge in $E^\prime$
must have both endpoints in $V^\prime$. Similarly for $E^{\prime \prime}$. 
But this implies that $G|_\sigma^c$ is disconnected, contrary to 
assumption.
So  $M$ is indecomposable of degree $-i$. By duality it also 
follows that we can construct the dual module, which has degree $2g-2+i$.

 It remains to show that there exist modules with degree in the range 
$0\le \mathrm{deg}(M)\le 2g-2$. Again by duality, it is enough to show the 
existence of modules with $0\le \mathrm{deg}(M)\le g$. Such modules are 
interesting, and their existence will be showed in the proposition below.
\end{proof}

The following result, which is also an analogue of a well-known fact from 
algebraic geometry, is needed to complete the proof of the previous proposition.

\begin{proposition}
Let $G$ be a graph of genus $g$. For any 
$0 \le d\le g$, there exists an indecomposable CM module $M$ 
with support $G$, locally of rank $1$, of degree $d$, multi-degree 
$\ge (0,0,\dots,0)$, and with a one-dimensional space of global sections,
i.e. $l(M) = 1$. 
\end{proposition}

\begin{proof}
If $d = 0$ we
can choose $M = k[G]$.  Otherwise $0 < d \leq g$.
Let $e_0 = \{v_0,w_0\}$ be en edge of $G$ whose removal gives a subgraph $H$ 
of genus $g-1$. By induction we may assume there
is a module $N$ on $H$ of degree $d-1$ and with the other properties
stated in the proposition.
 Make the module $M$ on $G$ such that $M_\sigma = N_\sigma$
when $\sigma \neq v_0, e_0$, but
$M_{e_0} = k\!\cdot \!1_{e_0}$, 
and $M_{v_0} = N_{v_0} \oplus k \!\cdot \!1_{v_0}$.
Define the $k$-linear maps $\varphi_{v e}$ of $M$ to be
general maps extending those of $N$, subject to the
commutativity constraint, i.e. they give an $S$-module
structure on $M$.
The degree of $M$ is $d$ and we see by Lemma \ref{RRLemCM} that
$M$ will be a CM module with $l(M)=l(N)=1$.
\end{proof}

\section{The Jacobian}\label{Jacobian}

We study the moduli space of isomorphism classes of squarefree CM $S$-modules,
with support on a   graph, of multi-degree $(0,0,\dots,0)$. 
This space is not as nice as the Jacobi variety of an algebraic curve. 
However, if we restrict to ``non-degenerate'' modules we can give the space an 
algebraic structure with dimension equal to the genus of the graph.

\begin{lemma} \label{RRLemTre}
Let $G$ be a tree, 
and $M$ a squarefree indecomposable module (possible not CM) of 
multi-degree $(0,0,\dots,0)$ with support $G$, locally of rank $1$.

For the vertices $v$ we may then choose generators $1_v$ of $M_v$ and similarly
$1_e$ for the edges $e$ such that whenever $v \in e$ then
$\varphi_{v e}$ sends $1_v$ to $1_e$.
In particular there is only one such CM module up to isomorphism, and
it will have a one-dimensional space of global sections.
\end{lemma}

\begin{proof}
We do induction on the number of vertices on $G$.
Let $v$ be a leaf of $G$, i.e. a vertex of degree one,
and $e = vw$ the incident edge.
Removing $v$ and $e$ we get a tree $H$. 
By induction it is enough to show that $\varphi_{w e}(1_w) = 1_e$ and 
$\varphi_{v e}(1_v) =1_e$. 
However, since both $\varphi_{w e}$ and $\varphi_{v e}$ are 
non-zero (if one of them is zero then $M$ is decomposable) we can first find 
a basis of $M_e$ such that $\varphi_{w e}=1$, and next find a basis for $M_v$ 
such that $\varphi_{v e}=1$. The same technique can be used to show the case 
where $G$ has only one edge. 

That $M$ has a one-dimensional space of global sections when $M$ is CM 
follows from Lemma \ref{RRLemCM}.
\end{proof} 

\begin{proposition}\label{RRProSykel}
Let $G$ be a cycle. 
The family of indecomposable CM modules over $k[G]$, locally of rank $1$, 
with multi-degree $(0,0,\dots,0)$ is parametrized by a union of 
$n$ $\mathbb{P}_k^1$'s, where we identify the same points from the different 
$\mathbb{P}_k^1$'s if they are not $(0,1)$ or $(1,0)$. In other words, it is
parametrized by a $\mathbb{P}_k^1$, with the points $(0,1)$ and $(1,0)$ 
$n$-doubled.
\end{proposition}

\begin{proof}
Suppose all the $\varphi_{ve}$ are nonzero when $v \in e$.
Let $e$ be an edge of $G$ and $H$ the line graph obtained by removing
$e$ from $G$.
 Then $M|_H$ has the $k$-linear maps as described in the lemma above. 
$M$ is therefore determined by the two $k$-linear maps $\varphi_{v e}$ and 
$\varphi_{w e}$, where $e=vw$. For any basis of $M_{e}$, these maps are 
determined by a pair $(s,t)\in k^2$. 
Since any basis change of $M_e$ is multiplication by a non-zero 
element of $k$, it follows that the maps are determined, up to isomorphism, 
by an element $(s,t)$ in $\mathbb{P}_k^1 \backslash \{ (1,0), (0,1) \}$.
Note that the latter is isomorphic to ${\mathbb A}_k^1 \backslash \{0\}$.

If two distinct $\varphi_{ve}$ are zero with $v$ incident to $e$
then $M$ is decomposable. So assume exactly one such $\varphi_{v_0e_0}$ is
zero. Then as in the case that $G$ is a tree that we may find generators
$1_v$ for the $M_v$ and $1_e$ for the $M_e$ such that all other
$\varphi_{ve}$ sends
$1_v$ to $1_e$ whenever $v \in e$. 
Hence for each $v_0 \in e_0$ there is exactly one isomorphism
class of such modules. If we have a cyclic order on the vertices, these
$2n$ pairs may be identified with $n$ copies of $(1,0)$ and $n$
copies of $(0,1)$.

In the cases of the first paragraph we will by Lemma \ref{RRLemCM} 
get CM modules with no global sections if $(s,t) \neq (1,1)$ in 
${\mathbb P}^1$ and with a one-dimension space of global sections, 
i.e. $l(M) = 1$ if $(s,t) = (1,1)$. 
In the cases of the second paragraph we get no global sections.


\end{proof}


\begin{proposition}
Let $G$ be a graph of genus $g\ge 1$. The family of indecomposable CM modules 
over $k[G]$, locally of rank $1$, with multi-degree $(0,0,\dots,0)$ and where
 all the maps $\varphi_{ve}$ are non-zero when $v \in e$,
is parametrized by $(\mathbb{A}^1_k\setminus \{0\})^g$.
\end{proposition}

\begin{proof}
We may choose a set $E^\prime$ of $g$ edges such that $(V, E \backslash E^\prime)$
is a tree. Given a module $M$ on $G$ of multidegree $(0, \ldots, 0)$ 
we may by Lemma \ref{RRLemTre} 
for the vertices $v$ choose generators $1_v$ of $M_v$ and
similarly generators $1_e$ for $M_e$ for the edges in $E \backslash E^\prime$
such that $\phi_{ve}(1_v) = 1_e$.  This choice of generators is unique
up to multiplication by a common scalar.

The module $M$ is now specified by for each edge $e = vw$ in $E^\prime$
choosing a pair of maps $\varphi_{ve}$ and $\varphi_{we}$. 
For a choice of generator $1_e$ of $M_e$ a pair of maps is given by 
a pair $(s,t)$ where $s,t \neq 0$. We may change the generator for $M_e$
by multiplying by a scalar, or change the set of generators $\{1_v\}$
by multiplying by a common scalar. This shows that such $M$ are
parametrized by  $(\mathbb{P}_k^1 \backslash \{ (1,0), (0,1) \})^g$
which is isomorphic to $({\mathbb A}_k^1 \backslash \{0\})$.


\end{proof}

\section{Gonality}\label{Go}
As we have seen, there is an analogy between line bundles on a given curve and 
squarefree Cohen--Macaulay modules with support on a given graph. 
We investigate this further by defining effective and special modules as an 
analogue of effective special divisors. The corresponding ``Brill--Noether'' 
theory for squarefree Cohen--Macaulay modules has some similarities and some 
differences from the classical theory.   

\begin{definition}
An indecomposable squarefree CM module $M$ supported on a graph $G$ 
is said to be 
effective if $M$ has a submodule isomorphic to $k[G]$, and said to be special 
if $M\subseteq \omega_{k[G]}$.
\end{definition}

A graph is CM if and only if it is connected, and $2$-CM if and only
 if it is $2$-connected. From Corollary \ref{2CM}, we know that the canonical
 module $\omega_{k[G]}$ is generated by the cycles $\widetilde{H}_1(G;k)$ if and
 only if $G$ is $2$-connected. It is therefore natural to study special modules
 of $2$-connected   graphs. 

Some of the following results are using that the field $k$ is infinite. 
We are therefore assuming that $\text{char}(k)=0$ for the rest of this section.
Recall that the cycles $Z_1(G)$ is the kernel of the boundary map from the
one-chains to the zero-chains
\[ C_1 = kE \pil kV = C_0. \]

\begin{lemma}\label{NinOmega}
Let $G$ be a two-connected   graph. Then $\omega_{k[G]}$ has a 
submodule $N\cong k[G]$ such that any squarefree module 
$N\subseteq M\subseteq \omega_{k[G]}$ is indecomposable, and therefore effective.
\end{lemma}

\begin{proof}
Let $u_N$ be a general element in $Z_1(G)$. Then $\text{Supp}(u_N)=G$. Let $N$ be the submodule of $\omega_{k[G]}$ generated by $u_N$ and suppose that $M$ is a module $N\subseteq M\subseteq \omega_{k[G]}$. Suppose that $M=M'\oplus M''$ is decomposable. Let $V'$ be the vector space spanned by the elements of $M'$ in degree $0$, and $V''$ the vector space spanned by 
the corresponding elements of $M''$. We will show that there exists a cycle $s\in Z_1(G)$ such that $s\not\in V'\oplus V''$. We construct $s$ as follows: Let $G^\prime := \mathrm{Supp} M^\prime$ and $G^{\prime \prime} :=\mathrm{Supp} M^{\prime \prime}$. We must have that $\dim G^\prime \cap G^{\prime \prime} < 1$. This is because $M_e \cong k$ for every edge $e$ in $G$, which means that each edge of $G$ must either be in $G^\prime$ or in $G^{\prime \prime}$. Since $G$ is connected, there has to be at least one vertex in this intersection. Let $v$ be any such vertex, and let $v^\prime$ be a vertex adjacent to $v$ in $G^\prime$ and $v^{\prime \prime}$ a vertex adjacent to $v$ in $G^{\prime \prime}$. Furthermore, since $G$ is two-connected, $G - \{v\}$ is connected and there is a simple path from $v^\prime$ to $v^{\prime \prime}$. Let $s\in Z_1(G)$ be a cycle corresponding to the simplicial circle obtained by connecting this path in $v$. This cycle is clearly not in $V^\prime \oplus V^{\prime \prime}$, since $s$ cannot be written as a sum of two cycles with
disjoint support. This means that $\dim_k V^{\prime}\oplus V^{\prime \prime} < \dim_k Z_1(G)$. 
Any decomposition of a module $M$ like this corresponds to a partitioning 
$E = E^{\prime}\cup E^{\prime \prime}$ of the edges in $G$, and in all cases, 
we will have as above 
that $\dim_k V^{\prime}\oplus V^{\prime \prime} < \dim_k Z_1(G)$. Since there are only
 finitely many ways to partitioning the edges into two sets, we can only find 
a decomposable $M$ between $N$ and $\omega_{k[G]}$ if $u_N$ lies in this finite 
union of subspaces of codimension $\ge 1$. But since $u_N$ is a general 
element, this is not the case.
\end{proof} 

\begin{definition}
Let $G$ be a $2$-connected graph of genus $g\ge 2$. A special and effective module $M$ satisfying $l(M)\ge r+1$ and $\mbox{deg}(M)=d$ is called a $g^r_d$.
 We define the gonality of a graph to be
$$\mbox{gon}(G)=\mbox{min }\left\{k\middle| G \text{ possesses a }g^1_k\right\}.$$ 
\end{definition}

For curves, it is well known that the gonality lies between $2$ and 
$\left\lfloor\frac{g+3}{2}\right\rfloor$. As we will see below, the same holds 
for a two-connected   graph.

\begin{lemma}\label{gonlegirth}
Let $G$ be a two-connected graph of genus $g\ge 2$. Then the gonality of $G$ 
is less then or equal to the girth 
(i.e. the length of the shortest nonzero cycle) of $G$.
\end{lemma}

\begin{proof}
From Lemma \ref{NinOmega}, we can choose $N\subseteq \omega_{k[G]}$ such that $N\cong k[G]$. Choose the cycle $u_N\in N$ and a cycle $c$ of length $w=\text{girth}(G)$. Then the submodule $M\subseteq \omega_{k[G]}$ generated by $u_N$ and $c$ is a $g^1_d$, where $d\le w$ since $\text{Supp}(u_N)\cap\text{Supp}(c)$ contains no more than $w$ vertices.
\end{proof}

\begin{lemma}\label{boundofgirth} Let $G$ be a two-connected graph of genus 
$g\ge 3$ and minimum valency $k\ge 3$. Then the girth of $G$ is
 $\le \left\lfloor\frac{g+3}{2}\right\rfloor$ except for three special cases:
\begin{itemize}
\item $K_{3,3}$,
\item The Petersen graph,
\item The Heawood graph.
\end{itemize}
\end{lemma}

\begin{proof}
According to \cite{Ca}, Theorem 11.11.3, we have that the number of vertices of
 $G$, denoted by $v$, is bounded by
$$2g-2\ge v \ge 1+k\frac{(k-1)^{\left\lfloor \frac{w-1}{2}\right\rfloor}-1}{k-2}.$$
Suppose that the girth $w\ge \left\lfloor \frac{g+5}{2}\right\rfloor$, then the 
inequality above has only finitely many solutions. The possible solutions are 
given as follows:
\begin{center}
\begin{tabular}{|c|c|c|}
\hline
$g$ & $v$ & $w$ \\
\hline
$3$ & $4$ & $4$ \\
\hline
$4$ & $4$--$6$ & $4$ \\
\hline
$6$ & $10$ & $5$ \\
\hline
$7$ & $10$--$12$ & $6$ \\
\hline
$8$ & $10$--$14$ & $6$ \\
\hline
$12$ & $22$ & $8$ \\
\hline
\end{tabular}
\end{center}


\noindent
However, the lower bound can be sharpened. If $w=4$ then $G$ must have a 
subgraph as in the picture below, as we now explain.
\begin{center}
\includegraphics[width=4cm]{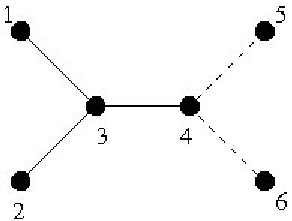}
\end{center}
Since $G$ has minimum valency $k\ge 3$, it is clear that it has a subgraph consisting of vertices $1$--$4$. But furthermore, since the vertex $4$ has valency $\ge 3$, and since the girth is $4$ there must exist two other vertices adjacent to $4$. Likewise, if $w=6$ and $G$ has minimum valency $k\ge 3$ then $G$ must contain a subgraph like the following figure.
\begin{center}
\includegraphics[width=5cm]{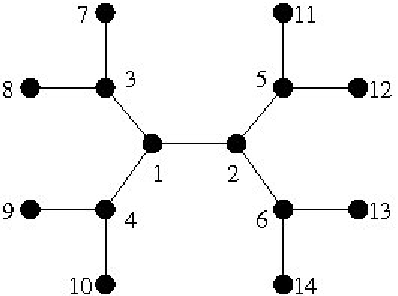}
\end{center}
The reason for this is the following. $G$ must contain a vertex, so $1$ is in $G$. Since $G$ has minimum valency $k\ge 3$, there has to be three vertices adjacent to $1$, call them $2,3$ and $4$. The vertex $2$ must also be adjacent to three vertices. But since the girth $w=6$ it is not adjacent to the vertex $2$ or $3$. Therefore it must be adjacent to two new vertices, call them $5$ and $6$. The same logic applies to the vertices $3,4,5$ and $6$, and it is clear that $G$ contains a subgraph as in the figure.  

For $w=8$, $G$ must also contain a subgraph as above. However, it is not so difficult to see, using the same argument as above, that it is impossible to construct a graph with $w=8$ and $k\ge 3$ with only $8$ more vertices. Therefore most of the solutions from the inequality above can not correspond to a graph. The only graphs satisfying the inequality in the Lemma that can exist have the following numerical data.

\begin{center}
\begin{tabular}{|c|c|c|}
\hline
$g$ & $v$ & $w$ \\
\hline
$4$ & $6$ & $4$ \\
\hline
$6$ & $10$ & $5$ \\
\hline
$8$ & $14$ & $6$ \\
\hline
\end{tabular}
\end{center}


There is exactly one graph for each of the cases above. The first is $K_{3,3}$, the second is the Petersen graph and the last is the Heawood graph.

\end{proof}

For any two-connected graph we would like to find a ``canonical'' model for it.
Let $G$ be a two-connected  graph of genus $g\ge 2$. If $v\in G$ is a vertex of
valency $2$ and $e_1=vw_1$ and $e_2=vw_2$ are the edges that meet $v$ and 
$e'=w_1w_2$ is not an edge in $G$. Then we construct a new graph 
$G'=(G-v)\cup\{e'=w_1w_2\}$. If we iterate this process, we end up with a 
graph we call $\widetilde{G}$. It is easy to verify that $\widetilde{G}$ has
the same topological structure as $G$. We will say that a graph $G$ is 
reduced if $G=\widetilde{G}$.

\begin{lemma}
Let $G$ be a two-connected graph of genus $g\ge 2$. 
Then $\emph{gon}(G)=\emph{gon}(\widetilde{G})$, and $\widetilde{G}$ is either
hyperelliptic, i.e. $2$-gonal, or has minimal valency $k\ge 3$.
\end{lemma}

\begin{proof}
Let $M$ be an effective module on $G$, and let $v$ be a vertex as described above. Then $k[G]_v\subseteq M_v\subseteq (\omega_{k[G]})_v$ but $k[G]_v\cong k \cong (\omega_{k[G]})_v$ since $v$ has valency $2$. This means that $M_v\cong k$, and it does not contribute to the degree of $M$. That means that for any effective module on $G$ we can construct an effective module $M'$ on $G'$ such that $M'$ has the same numerical invariants as $M$. Iterating this process we get an effective module $\widetilde{M}$ on $\widetilde{G}$ with the same numerical invariants as $M$. It is also easy to see that any effective module $N$ on $\widetilde{G}$ can be obtained from an effective module coming from $G$. This proves the first part.

For the other part, we notice that either $\widetilde{G}$ has minimal valency 
$k\ge 3$ or there exists a vertex $v$ of valency $2$ such that if $e_1=vw_1$ 
and $e_2=vw_2$ are the edges that meet $v$, then $e'=w_1w_2$ is an edge in 
$\widetilde{G}$. That means that there exists a cycle $s\in Z^1(G)$ such that 
$\text{Supp}(s)=e_1\cup e_2\cup e'$. Since $G$ is $2$-connected we can find a 
submodule $N\subseteq\omega_{k[G]}$, isomorphic to $k[G]$, such that any module 
in between is effective. Let $u_N$ be the generator of $N$, then $M=(u_N,s)$ is
 a $g^1_d$ where $d\le 3$ since $\text{Supp}(u_N)\cap\text{Supp}(s)$ contains 
$3$ vertices. Furthermore, $d\le 2$ since $v$ has valency $2$, and $d\ge 2$ 
sine $G$ is two-connected.
Hence $\widetilde{G}$ has gonality $2$, that is $\widetilde{G}$ is
hyperelliptic.
\end{proof}

\begin{theorem}\label{gonbound}
Let $G$ be a two-connected graph of genus $g\ge 2$. 
Then $2\le \mathrm{gon}(G)\le \left\lfloor\frac{g+3}{2}\right\rfloor$.
\end{theorem}
\begin{proof}
That the gonality cannot be $1$ is clear since $G$ is two-connected.
By the previous two Lemmata, it is enough to show that there exists a 
$g^1_{\left\lfloor\frac{g+3}{2}\right\rfloor}$ on the three special graphs above. 
The graph $K_{3,3}$ can be drawn like this:
\begin{center}
\includegraphics[width=4cm]{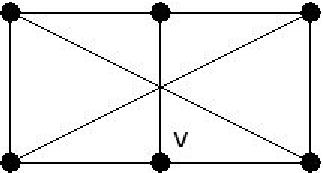}
\end{center}
We notice the two subgraphs
\begin{center}
$G_1$: \includegraphics[width=3cm]{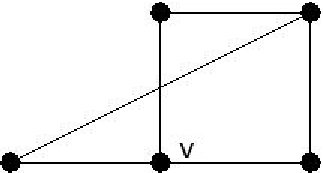} \quad \quad and \quad \quad
$G_2$: \includegraphics[width=3cm]{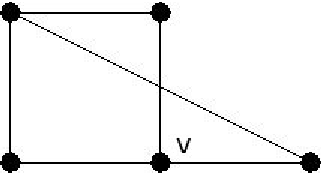}
\end{center}
Let $s_1\in Z^1(G)$ be a cycle with support on $G_1$. We claim that there exists a cycle $s_2\in Z^1(G)$ with support on $G_2$ such that the module $M$ 
generated by $s_1$ and $s_2$ is a $g_3^1$. First of all we notice that for any 
$s_2$ with support on $G_2$ we have that $M$ is effective since 
$\text{Supp}(s_1,s_2)=G$, and indecomposable since $\text{dim } G_1\cap G_2=1$. 
Next, we show that we can choose $s_2$ such that the restrictions of $s_1$ and 
$s_2$ are equal in the vertex $v$. This follows if we show that we can choose 
$s_2$ such that the restrictions of $s_1$ and $s_2$ are equal in the three 
edges of $v$, and
this is so since the two four-cycles of $G_2$ are linearly independent
at the vertex $v$.
We now have
$$\text{deg}(M)=\sum_{\text{vertices} v} 
(\text{dim}_k\left(\text{Span}\{s_1|_{v},s_2|_{v}\}
\right)-1)=3.$$

The Petersen graph is given as:
\begin{center}
\includegraphics[width=4cm]{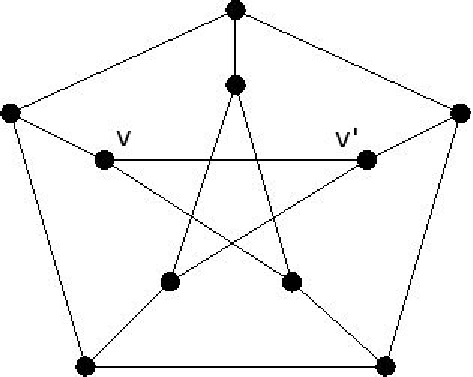}
\end{center}
We notice the two subgraphs
\begin{center}
$G_1$: \includegraphics[width=3cm]{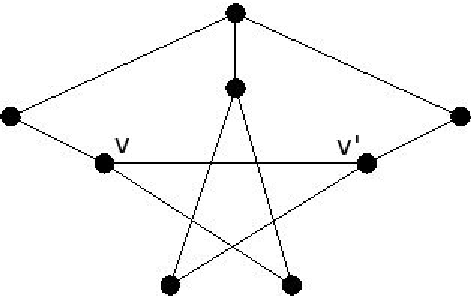} \quad \quad and \quad \quad
$G_2$: \includegraphics[width=3cm]{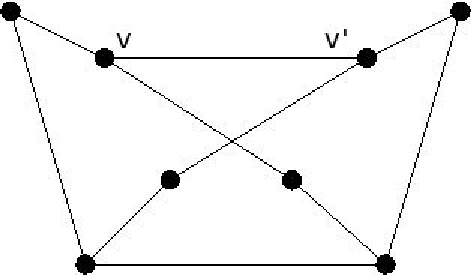}
\end{center}
Let $s_1$ be a cycle with support on $G_1$. As above, we claim that there 
exists a cycle $s_2$ with support on $G_2$ such that the module $M$ generated 
by $s_1$ and $s_2$ is a $g_4^1$. For any $s_2$ with support on $G_2$ we have 
that $M$ is effective and indecomposable for the same reasons as above. Next, 
we show that we can choose $s_2$ such that the restrictions of $s_1$ and $s_2$
 are equal in the vertices $v$ and $v'$. As above, this follows if we can 
choose $s_2$ such that the restrictions of $s_1$ and $s_2$ are equal on any two
 of the edges of $v$ and any two of the edges of $v'$. Again, this follows 
since  $G_2$ has three linearly independent cycles, $v$ and $v^\prime$
each has valency $3$, and cycles are independent in these vertices.
Similarly as above, we now get 
that $\text{deg}(M)=4$.

The Heawood graph is given as:
\begin{center}
\includegraphics[width=4cm]{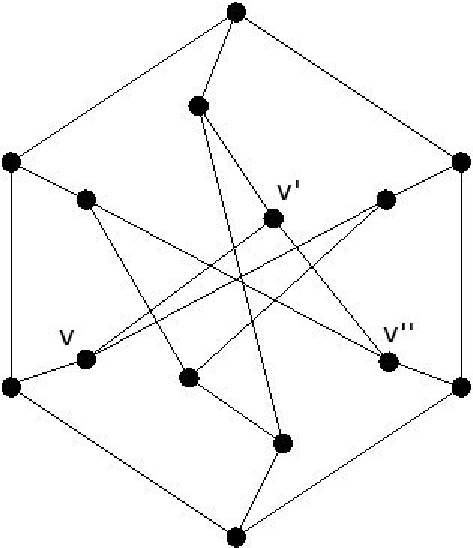}
\end{center}
We notice the two subgraphs
\begin{center}
$G_1$: \includegraphics[width=3cm]{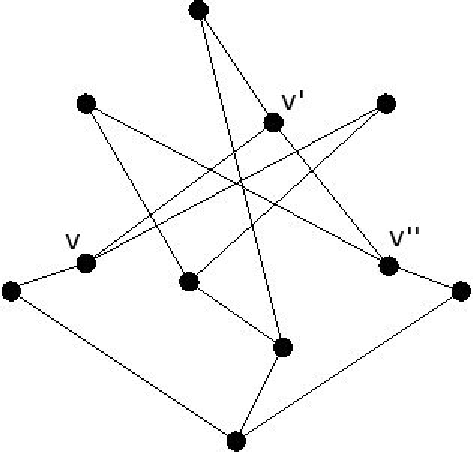} \quad \quad and \quad \quad
$G_2$: \includegraphics[width=3cm]{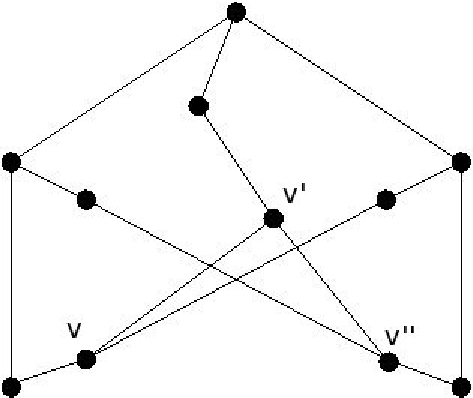}
\end{center}
Let $s_1\in\widetilde{H}_1(G;k)$ be a cycle with support on $G_1$. As above, we claim that there exists a cycle $s_2\in\widetilde{H}_1(G;k)$ with support on $G_2$ such that the module $M$ generated by $s_1$ and $s_2$ is a $g_5^1$. For any $s_2$ with support on $G_2$ we have that $M$ is effective and indecomposable for the same reasons as above. Next, we show that we can choose $s_2$ such that the restrictions of $s_1$ and $s_2$ are equal in the vertices $v,v'$ and $v''$. As above, this follows if we can choose $s_2$ such that the restrictions of $s_1$ and $s_2$ are equal on any two of the edges of $v,v'$ and $v''$. Again, this follows since $G_2$ has four linearly independent cycles and the three vertices
$v, v^\prime$, and $v^{\prime \prime}$ each have valency $3$.
Similarly as above, we now get that $\text{deg}(M)=5$.

\end{proof}
One may ask which gonalities that occur in the bound. It is well-known that in
 the moduli space of curves of genus $g$, there exist curves of gonality $k$ 
for all $2\le k\le \left\lfloor\frac{g+3}{2}\right\rfloor$. However, for the 
finite set of reduced graphs of genus $g$ this is not the case. For
 large enough $g$ the upper bound is not sharp. This can be seen by combining
 Lemma \ref{gonlegirth} with the inequality that bounds the number of vertices
 in the beginning of the proof of Lemma \ref{boundofgirth}.



\section{Clifford's theorem} \label{Cl}

We shall now prove the analog of Clifford's theorem. This is a theorem
concerning vector spaces of cycles on a graph, and so may be formulated
without the notion of squarefree modules. In algebraic geometry 
Clifford's theorem concerns linear systems on smooth
projective curves. Graphs
are in some sense singular, and in order for Clifford's theorem to 
hold we must have an extra assumption on our space of cycles. It must
fulfill the criterion (to be defined) of being {\it locally indecomposable} 
at each vertex of the graph.

\subsection{Definitions}
We recall some notions from graph theory which will be necessary, and
make several definitions needed in the statements and proofs of 
Clifford's theorem.

Let $G$ be a connected graph. A vertex $P$ of $G$ is called a 
{\it cut point} of $G$ if the removal of $P$ disconnects $G$. The
maximal subgraphs of $G$ which have no cut points are called the blocks
of $G$. A block consists either just of an edge, or it is a $2$-connected
graph. Note that each block of $G$ must contain at least one of the cut 
points of 
$G$. Those blocks which contain only one cut point will be called
{\it leaf blocks} of $G$. The number of leaf blocks of $G$ is
denoted by $l(G)$. If $Q$ is a cut point of $G$, then the removal of
$Q$ and the incident edges disconnects $G$ into components 
$\hat{G}_1, \hat{G}_2, \ldots, \hat{G}_r$. Adding to $\hat{G_i}$ 
the edges in $G$ between $\hat{G}_i$ and $Q$, we get subgraphs
$G_1, G_2, \ldots, G_r$ such that $G$ is their union and 
any pair intersects only in the cut point $Q$. 
These are the {\it cut components}
of $G$ at the cut point $Q$. Note that $Q$ will not be a
cut point of any of the $G_i$. 

  Let $G=(V,E)$ where $V$ are the vertices of $G$ and $E$ its edges. 
Let $C_1 = kE$ be the one-chains of $G$. If $c$ is a chain, we define
the support $\supp c$ to be the induced subgraph on the set of edges in $c$.
The cycles $Z_1(G)$ is the kernel of the boundary map from the one-chains
to the zero-chains
\[ C_1 = kE \pil kV = C_0. \]

\begin{lemma} Let $Q$ be a cutpoint of $G$ and $G_1, \ldots, G_r$ its
cut components at $Q$. If $s$ is a cycle on $G$, then 
$s = s_1 + s_2 + \cdots + s_r$ where $s_i$ is a cycle on $G_i$.
Hence there is an isomorphism
\[  Z_1(G) \iso \bigoplus\limits_{i = 1}^r Z_1(G_i). \]
\end{lemma}

\begin{proof} Write $s = s_1 + \cdots + s_r$ where $s_i$ is a chain in 
$G_i$. Then
\[ 0 = \vardel (s) = \vardel(s_1) + \cdots + \vardel(s_r). \]
Now $\vardel(s_i)$ is a zero chain in $G_i$. Every point $P$ involved
in $\vardel(s_i)$ must cancel against some point in $\vardel(s_j)$ for 
some $j \neq i$. But since $G_i$ and $G_j$ only intersect in $Q$
we must have $\vardel(s_i) = \alpha_i Q$ for a constant $\alpha_i$. 
Since $\vardel^2(s_i) = 0$, we get each $\alpha_i = 0$, and so 
the $s_i$ are cycles in $G_i$. 
\end{proof}

A subspace  $\Gamma \sus Z_1(G)$ of the cycles
will be called a {\it cycle system}. The support $\supp \Gamma$
is the union of the supports of all cycles in $\Gamma$. 
We shall assume the support of $\Gamma$ is $G$. Then each block of 
$G$ must be two-connected.
In the case of the lemma above, 
let $\Gamma_i$ be the image of $\Gamma$ in $Z_1(G_i)$. We then get an
injection 
\begin{equation} \Gamma \inpil \bigoplus\limits_{i = 1}^r \Gamma_i. 
\label{CliffLigGamma} \end{equation}

\medskip

For each vertex $P \in V$, let $E_P$ be the set of edges incident to $P$
and let $C_{1P} = kE_P$ be the vector space with these edges as basis.
There is a natural map
\[ C_{1P} = kE_P \pil kP \]
sending each edge $e \mapsto P$, and let $Z_{1P}(G)$ be the kernel of this map,
the local cycles at $P$. There is a natural map $C_1 \pil C_{1P}$ and
$Z_1(G)$ maps to $Z_{1P}(G)$. 

 If $\Gamma$ is a cycle system, we get for each vertex $P$ a map
$\Gamma \pil Z_{1P}(G)$, and denote by $\Gamma_P$ the image of this.
We define the degrees
\begin{eqnarray*}  d(\Gamma,P) & = &  \dim_k \Gamma_P - 1 \\
  d(\Gamma) & = & \sum_{P \in V} d(\Gamma,P).
\end{eqnarray*}

Now $\Gamma_P \sus Z_{1P}(G)$ and this is again a subspace of $C_{1P} = kE_P$.
The latter is a vector space with a natural basis, the edges incident to
$P$. Let us now for a moment consider this situation.
Let $V$ be a vector space over the field $k$ with a basis $B = \{b_1,
\ldots, b_n\}$. Let $\pi = B_1 \cup B_2 \cup \cdots \cup B_r$ be 
a partition of $B$ into nonempty parts. 
A subspace $U \sus V$ is {\it decomposable} with respect
to $\pi$ if 
\[ U \iso (U \cap kB_1) \oplus (U \cap kB_2) 
\oplus \cdots \oplus (U \cap kB_r). \]
We call $\Gamma$ {\it locally decomposable at the vertex $P$} is 
there is a partition of $E_P$ such that $\Gamma_P$ is decomposable
with respect to this partition;
otherwise $\Gamma$ is  {\it locally indecomposable at the vertex $P$}.
We say that $\Gamma$ is  {\it locally indecomposable}
if it is locally indecomposable at all its vertices.
We say $\Gamma \sus kE$ is {\it decomposable} if there is a partition of 
the set of edges such that $\Gamma$ is decomposable with respect to 
this partition;
otherwise $\Gamma$ is {\it indecomposable}.

\medskip
Let $P$ be a vertex where $\Gamma$ is decomposable.
Let $E_P = E_{P1} \cup E_{P2}$ be a 
partition with respect to which $\Gamma_P$ is decomposable.
We may now make a new graph $G^\prime$ by replacing the 
vertex $P$ by two vertices $P_1$ and $P_2$  and let the edges in
$E_{P1}$ be incident to $P_1$ and those in $E_{P2}$ be incident to $P_2$. 
There is an exact sequence
\[ 0 \pil K_P \pil \Gamma \pil \Gamma_P \pil 0 \]
where $K_P$ is the kernel. Note that 
\[ \Gamma_P = (\Gamma_P \cap kE_{P_1}) \oplus (\Gamma_P \cap kE_{P_2}). \]
Define the kernels
\begin{eqnarray*}
\Gamma_1 & = & \ker (\Gamma \pil (\Gamma_P \cap k E_{P_2})) \\
\Gamma_2 & = & \ker (\Gamma \pil (\Gamma_P \cap k E_{P_1})).
\end{eqnarray*} 
For $i = 1,2$ chose splittings of the surjections 
\[ \Gamma_i \overset{\overset{\sigma_i}\leftarrow}{\lpil} 
(\Gamma_P \cap k E_{Pi}). \]
We then have a decomposition
\begin{equation}
\Gamma \iso K_P \oplus \im \sigma_1 \oplus \im \sigma_2. 
\label{CliffLigGammasplitt}
\end{equation}
There is an injection $Z_1(G^\prime) \inpil Z_1(G)$ and we see that the image
of this map contains all the summands to the right of the above
isomorphism. We then get a cycle system $\Gamma^\prime$ on $G^\prime$
mapping isomorphically onto $\Gamma$. This subspace $\Gamma^\prime$
does not depend on the choice of splittings.

If $\Gamma^\prime$ is again locally decomposable at $P_1$ or $P_2$
we may continue the process. In the end we get a graph $\tilde{G}$
and a cycle system $\tilde{\Gamma}$ which is locally indecomposable
at all vertices $Q$ of $\tilde{G}$ mapping to $P$. 
We call $(\tilde{G}, \tilde{\Gamma})$ the {\it resolution} of 
the pair $(G, \Gamma)$ at the vertex $P$.

\subsection{Versions of Clifford's theorem}

\medskip
Clifford's theorem in algebraic geometry relates the dimension and
degree of a linear system
on a projective curve.

\begin{theorem}[Clifford's theorem] Let $D$ be an effective special divisor  
on a smooth projective curve $C$ over $\Spec k$. 
Let $\gL(D)$ be the associated line bundle, and
$\Gamma (\gL(D))$ its global sections.
Then $\deg D + 2 \geq 2 \dim_k \Gamma (\gL(D))$. 
\end{theorem}



Here is our version for cycle systems on graphs.

\begin{theorem}[Clifford's theorem] \label{CliffTheToshg}
Let $G$ be a two-connected graph and $\Gamma$ a cycle system with support
$G$. If $\Gamma$ is locally indecomposable, then 
\[ d(\Gamma) + 2 \geq 2 \dim_k \Gamma. \]
\end{theorem}

\begin{example}
The following example shows that even if $\Gamma$ is indecomposable,
it is necessary that $\Gamma$ be locally indecomposable for Clifford's 
theorem to hold. Let $G$ be the two-connected graph
\begin{center}
\includegraphics[width=8cm]{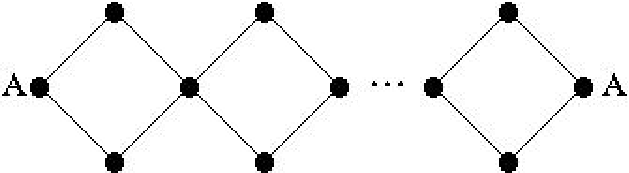}
\end{center}
\noindent where the endpoints $A$ are identified. 
Let $c_1, c_2, \ldots, c_r$ be the four-cycles in the graph, and
let $\Gamma$ be generated by the cycles
\[ c_1 - c_2, c_2 - c_3, \ldots, c_{r-1} - c_r. \]
Then $\Gamma$ is indecomposable. But $\Gamma_P$ is locally decomposable at
each vertex $P$. We have $d(\Gamma) = r$ and $\dim_k \Gamma = r-1$
so the inequality $d(\Gamma) + 2 \geq 2 \dim_k \Gamma$ far from holds.
\end{example}

\begin{corollary} \label{CliffCorTrivalens}
Let $G$ be a two-connected graph where every vertex has degree at most
three, and $\Gamma$ a cycle system with support $G$. Then
\[ d(\Gamma) + 2 \geq 2 \dim_k \Gamma. \]
\end{corollary}

\begin{proof}
In this case $\Gamma_P$ cannot be decomposable for any vertex $P$, so 
$\Gamma$ is locally indecomposable.
\end{proof}

We will prove Theorem \ref{CliffTheToshg} through an inductive
argument. It will then be necessary to have versions for the
cases when the cycle systems may be locally decomposable.

If $G$ is a graph and $\Gamma$ a cycle
system on $G$ we denote by $q(\Gamma)$ the number of vertices of $G$
at which $\Gamma$ is locally decomposable. 
When $G$ has connectivity one a leaf block of $G$ will be called {\it regular}
if all its vertices are locally indecomposable,
except perhaps the vertex to which it is attached to $G$ (the cutpoint). 
Let $l_{reg}(\Gamma)$ be the number of regular leaf blocks of $G$.

\begin{proposition} \label{CliffProConOne} Let $G$ be a graph of connectivity 
one, and
$\Gamma$ a cycle system with support $G$. Then
\[ d(\Gamma) + q(\Gamma) + l_{reg}(\Gamma) \geq 2 \dim_k \Gamma. \]
\end{proposition}

\begin{proposition} \label{CliffProConTwo} Let $G$ be
a two-connected graph, and $\Gamma$ a cycle system with support $G$.
Suppose the number of locally decomposable vertices is $q(\Gamma)
\geq 1$. Then
\[ d(\Gamma) + q(\Gamma) \geq 2 \dim_k \Gamma. \]
\end{proposition}

In order to formulate Clifford's theorem for modules, we define a $k[G]$-module
$M$ to be indecomposable at a vertex $i$ if the submodule $M_{\geq i}$
of $M$, 
consisting of the multigraded pieces of $M$ indexed by multidegrees $\geq i$, 
is an indecomposable module. The module $M$ is locally indecomposable
if it is indecomposable at each vertex $i$ of $G$. 
Formulated in terms of modules we get.
\begin{theorem}[Clifford's theorem] 
Let $G$ be a two-connected graph and $M$ a special effective 
$k[G]$-module. If $M$ is locally indecomposable, then 
\[ \deg(M) + 2 \geq 2 \, l(M). \]
\end{theorem}

\begin{proof}
Let $N$ be the submodule of $M$ generated by the sections $M_{\emptyset}$.
Then $N$ is also special and effective. Let $q$ be the number
of vertices at which $N$ is locally decomposable.
Then $\deg(M) \geq \deg(N) + q$. If $q = 0$ the statement follows by
Theorem \ref{CliffTheToshg}. If $q \geq 1$ we have by Proposition
\ref{CliffProConTwo}
\[ \deg(M) \geq \deg(N) + q \geq 2 \, l(N) = 2 \, l(M), \]
which implies the inequality of the theorem.
\end{proof}

In the same way as we got Corollary \ref{CliffCorTrivalens} we
get the following.

\begin{corollary} Let $G$ be a two-connected graph where every vertex
has degree at most three, and $M$ a special
effective $k[G]$-module. Then
\[ \deg(M) + 2 \geq 2 \, l(M). \]
\end{corollary}
\subsection{The proofs}

We will prove Theorem \ref{CliffTheToshg} and Propositions 
\ref{CliffProConOne} and \ref{CliffProConTwo} by induction on the number of
locally decomposable vertices and the number of edges of $G$. The
argument will be simultaneous induction in the sense that when we
prove one of the statements we assume that all three statements
hold if the number of edges is strictly smaller or the number
of edges is the same, but the number of locally decomposable
vertices is strictly smaller.

\begin{proof}[Proof of Proposition \ref{CliffProConOne}.] Let $Q$
be a cut point of $G$. If $F$ is a cut component of $G$ at $Q$,
then $\Gamma$ induces a cycle system $\Phi$ on $F$.
For all points $P$ in $F$ distinct from $Q$ we will have $\Phi_P
= \Gamma_P$. Note that the point $Q$ on $F$ will not be a cut point of $F$.

If $F$ has connectivity one, either $Q$ will be in a
nonleaf block of $F$, denote the set of such $F$'s by $NL^1$,
or $Q$ will be in a leaf block of $F$. This leaf block may
either be regular, denote the set of such $F$'s by $L^1_{reg}$,
or non-regular, denote the set of such $F$'s by $L^1_{nreg}$.
Let $C^1$ be the union of these three sets, the set of
$F$'s of connectivity one.

Let now  $F$ be two-connected. Either $\Phi$ is locally
indecomposable, denote the set of such $F$'s by $B^2_{reg}$,
or the only locally decomposable vertex for $\Phi$ is
$Q$, denote the set of such $F$'s by $B^2_Q$, or
$\Phi$ contains a vertex distinct from $Q$ which is
locally decomposable, denote the set of such $F$'s by $B^2_{nreg}$.
Let $C^2$ be the union of these three sets, the set of
two-connected $F$'s. Finally let $C = C^1 \cup C^2$ be
the set of all the components $F$.

By induction, when $F$ has connectivity one, we have 
\begin{equation} \label{CliffLig1}
d(\Phi) + q(\Phi) + l_{reg}(\Phi) \geq 2 \dim_k \Phi,
\end{equation}
and
when $F$ is two-connected we have
\begin{alignat}{2} \label{CliffLig2reg}
d(\Phi) + 2 & \, \geq \, & 2 \dim_k \Phi, \quad & (F,\Phi) \in B^2_{reg}\\
d(\Phi) + q(\Phi) & \, \geq \, & 2 \dim_k \Phi, \quad & 
(F,\Phi) \in B^2_Q \cup B^2_{nreg}. \label{CliffLig2nreg}
\end{alignat}

If $F$ has connectivity one and is in $NL^1$ or $L^1_{nreg}$,
the regular leaf blocks of $F$ will still be regular leaf blocks
in $G$, but when $F$ is in $L^1_{reg}$, the attaching regular
leaf block will cease to be a leaf block. The two-connected
$F$'s which are in $B^2_{reg}$ and $B^2_Q$ will become regular
leaf blocks of $G$, but not the ones in $B^2_{nreg}$.
We therefore have
\begin{eqnarray} \notag  l_{reg}(\Gamma) &\,  = \, & \sum_{(F,\Phi) \in NL^1}
l_{reg}(\Phi) + \sum_{(F,\Phi) \in L^1_{nreg}} l_{reg}(\Phi) \\ 
\label{CliffLigL} & + & \sum_{(F,\Phi) \in L^1_{reg}} (l_{reg}(\Phi) - 1) + |B^2_{reg}| + |B^2_Q| \\
& \, = \, & \sum_{(F,\Phi) \in C^1} l_{reg}(\Phi) + |B^2_{reg}| + |B^2_Q|
- |L^1_{reg}|. \notag
\end{eqnarray}

Concerning the locally decomposable vertices with respect to $\Phi$,
when $F$ is in $NL^1$ or $L^1_{nreg}$ we may lose a decomposable vertex
when going from $\Phi$ to $\Gamma$, since $Q$ may be changing its status 
from locally decomposable to locally indecomposable.
Similarly for $F$ in $B^2_Q \cup B^2_{nreg}$.
Hence
\begin{eqnarray} \notag
q(\Gamma) & \geq & \sum_{(F,\Phi) \in NL^1 \cup L^1_{nreg}} (q(\Phi)
-1) + \sum_{(F,\Phi) \in L^1_{reg}} q(\Phi) \\ \label{CliffLigQ}  & + & 
\sum_{(F,\Phi) \in B^2_Q  \cup B^2_{nreg}} (q(\Phi) - 1) + \epsilon \\
& = & \sum_{(F,\Phi) \in C} q(\Phi) -|NL^1| - |L^1_{nreg}| \notag
- |B^2_Q| - |B^2_{nreg}| + \epsilon, 
\end{eqnarray}
where $\epsilon = 1$ if $\Gamma$ is locally decomposable at $Q$, and $\epsilon
= 0$ if $\Gamma$ is locally indecomposable at $Q$.

By (\ref{CliffLigGamma}) there is a commutative diagram
\begin{equation} \label{CliffCD}
\begin{CD} \Gamma @>>> \oplus_{(F,\Phi) \in C} \Phi @>>> T \\
@VVV @VVV @VVV \\
\Gamma_Q @>>> \oplus_{(F, \Phi) \in C} \Phi_Q @>>> T_Q
\end{CD}
\end{equation}
where $T$ and $T_Q$ are the cokernels of the injections on the left.
Since each $\Phi \pil \Phi_{Q}$ is surjective, there is 
a surjection $T \pil T_Q$. 

From this we obtain

\begin{eqnarray}
d(\Gamma, Q) + 1 & = & \sum_{(F,\Phi) \in C} (d(\Phi,Q) + 1)
- \dim T_Q \\
d(\Gamma,Q) & = & \sum_{(F,\Phi) \in C} d(\Phi,Q)
+|C| - 1 - \dim_k T_Q.
\end{eqnarray}
At all points $P$ of $F$ distinct from $Q$ we have $d(\Gamma,P) =
d(\Phi,P)$. Therefore
\begin{equation} \label{CliffLigD}
 d(\Gamma) = \sum_{(F,\Phi) \in C} d(\Phi) + |C|
-1 - \dim_k T_Q.
\end{equation}
Now we add the equations obtained by induction (\ref{CliffLig1}),
(\ref{CliffLig2reg}), and (\ref{CliffLig2nreg}) and get
\[ \sum_{(F,\Phi) \in C} d(\Phi) + \sum_{(F,\Phi) \in C} q(\Phi)
+ 2|B^2_{reg}| + \sum_{(F,\Phi) \in C^1} l_{reg}(\Phi) \geq
2 \sum_{(F,\Phi) \in C} \dim_k \Phi. \]

Equations (\ref{CliffLigD}), (\ref{CliffLigQ}), and
(\ref{CliffLigL}) give expressions for 
$d(\Gamma), q(\Gamma)$ and $l_{reg}(\Gamma)$ respectively.
Inserting these in the above equation and taking
into account the upper row in the commutative diagram
(\ref{CliffCD}), we obtain
\begin{eqnarray*}
& & d(\Gamma) - |C| + 1 + \dim_k T_Q \\
& + & q(\Gamma) + |NL^1| + |L^1_{nreg}| + |B^2_Q| + |B^2_{nreg}|
- \epsilon \\
& + & 2|B^2_{reg}| + l_{reg}(\Gamma) + |L^1_{reg}| - |B^2_Q| - |B^2_{reg}| \\
& \geq & 2 \dim_k \Gamma + 2 \dim_k T.
\end{eqnarray*}

The terms containing the cardinalities of sets of cut components are
\[ - |C| + |NL^1| + |L^1_{nreg}| + |L^1_{reg}| + |B^2_{reg}|
+ |B^2_{nreg}| \leq 0. \]

Hence we get
\[ d(\Gamma) + l_{reg}(\Gamma) + q(\Gamma) \geq
2 \dim_k \Gamma + (2 \dim_k T - \dim T_Q + \epsilon - 1). \]
If $\dim T \geq 1$, the expession in the paranthesis on the right
is $\geq 0$ since $\dim_k T \geq \dim_k T_Q$.
If $\dim_k T = 0$ then $\epsilon = 1$ and so the paranthesis is also
$\geq 0$. Hence we obtain
\[ d(\Gamma) + l_{reg}(\Gamma) + q(\Gamma) \geq 2 \dim_k \Gamma. \]
\end{proof}

\begin{proof}[Proof of Proposition \ref{CliffProConTwo}.]
Let $P$ be a locally decomposable vertex on $G$ for the cycle system
$\Gamma$ and let $(\tilde{G}, \tilde{\Gamma})$ be the resolution of 
$(G, \Gamma)$ at the vertex $P$. Since $G$ is two-connected,
$\tilde{G}$ will be connected.

If $\tilde{G}$ is two-connected we will have $q(\tilde{\Gamma}) = q(\Gamma) - 1$
and $d(\Gamma) \geq d(\tilde{\Gamma}) + 1$. If $q(\tilde{\Gamma}) \geq 1$
we have by induction 
\[ d(\tilde{\Gamma}) + q(\tilde{\Gamma}) \geq 2 \dim_k \tilde{\Gamma}, \] 
and so we obtain
\begin{eqnarray*}
& &  d(\Gamma) + q(\Gamma) \geq d(\tilde{\Gamma}) + q(\tilde{\Gamma}) + 2 \\
& \geq & 2 \dim_k \tilde{\Gamma} + 2 = 2 \dim_k \Gamma + 2 \geq 2 \dim_k
\Gamma.
\end{eqnarray*}

If $q(\tilde{\Gamma}) = 0$ we have by induction
\[ d(\tilde{\Gamma}) + 2 \geq 2 \dim \tilde{\Gamma} \]
so 
\[ d(\Gamma) + q(\Gamma) = d(\Gamma) + 1 \geq d(\tilde{\Gamma}) + 2 
\geq 2 \dim_k \Gamma. \]

\medskip
Suppose now $\tilde{G}$ has connectivity one.
Then 
\[ d(\tilde{\Gamma}) + q(\tilde{\Gamma}) + l_{reg}(\tilde{\Gamma})
\geq 2 \dim_k \tilde{\Gamma}. \]
Let $D$ be the set of vertices in $\tilde{G}$ which map to 
$P$ in $G$. Then each leaf block of $\tilde{G}$ must contain an element of
$D$ which is not the vertex at which the leaf block is attached
to $\tilde{G}$ (the cut point). Otherwise the leaf
block maps to a leaf block of $G$, which is impossible since $G$ is 
two-connected. Thus $|D| \geq l(\tilde{G})$. 

Now we have $d(\Gamma) = d(\tilde{\Gamma}) + |D| - 1$. Also
$q(\Gamma) = q(\tilde{\Gamma}) + 1$. We therefore get
\begin{eqnarray*}
d(\Gamma) + q(\Gamma) = d(\tilde{\Gamma}) + q(\tilde{\Gamma}) + |D| 
& \geq & d(\tilde{\Gamma}) + q(\tilde{\Gamma}) + l(\tilde{G}) \\
& \geq & d(\tilde{\Gamma}) + q(\tilde{\Gamma}) + l_{reg}(\tilde{\Gamma})\\
& \geq & 2 \dim_k \tilde{\Gamma} = 2 \dim_k \Gamma.
\end{eqnarray*}
\end{proof}

Finally we are able to prove Clifford's theorem.

\begin{proof}[Proof of Theorem \ref{CliffTheToshg}.]
Choose a cycle $s$ in $\Gamma$, and an edge $e$ on $s$.
Let $\Gamma_0$ be the subspace of $\Gamma$ consisting of cycles not
containing the edge $e$. Then $\Gamma = \Gamma_0 \oplus ks$. Let $G_0$
be the support of $\Gamma_0$, so $G_0$ has less edges than $G$. 
On each connected component $F$ of $G_0$ the cycle system $\Gamma_0$ induces
a cycle system $\Phi$. 
Denote by $C^1$ the set of components of $G_0$ of connectivity one,
and by $C^2$ the set of components of $G_0$ which are two-connected.
We may decompose $C^2 = C^2_{reg} \cup C^2_{nreg}$ where
$C^2_{reg}$ are the components $F$ of $G_0$ where the induced cycle system
$\Phi$ is locally indecomposable, and $C^2_{nreg}$ are those
where $\Phi$ has some vertex which is locally decomposable.

Now we have 
\begin{eqnarray}
d(\Phi) + l_{reg}(\Phi) + q(\Phi) & \geq & 2 \dim_k \Phi, \quad (F,\Phi) \in C^1 
\notag \\ \label{CliffLigInduk}
d(\Phi) + 2 & \geq & 2 \dim_k \Phi, \quad (F,\Phi) \in C^2_{reg} \\ \notag
d(\Phi) + q (\Phi) & \geq & 2 \dim_k \Phi, \quad (F,\Phi) \in C^2_{nreg}.
\end{eqnarray}


Let $T$ be the set of vertices $P$ of $G_0$ where 
$\dim_k \Gamma_P > \dim_k \Gamma_{0P}$.
We may decompose $T = T^1 \cup T^2$ where $T^1$ are the vertices in $T$
contained in components of connectivity one and $T^2$ the vertices
which are in two-connected components. 
Furthermore $T^1 = T^1_{reg} \cup T^1_{nreg}$ where $T^1_{reg}$ are
the points which are in regular leaf blocks, and are not cutpoints
in the component they belong to, and $T^1_{nreg}$ are the rest
of the points in $T^1$. 
Also $T^2 =  T^2_{reg} \cup T^2_{nreg}$ where $T^2_{reg}$ are the points
of $T$ which are on components in $C^2_{reg}$ and $T^2_{nreg}$
are the points of $T$ which are on components of $C^2_{nreg}$.


As we now explain the following three inequalities hold:
\begin{eqnarray}
\notag |T^1_{reg}| &\geq & \sum_{(F,\Phi) \in C^1}  l_{reg}(\Phi),  
\notag \\ \label{CliffLigUlik} 
|T^{1}_{nreg} \cup T^2_{nreg}| & \geq & \sum_{(F,\Phi) \in C} q(\Phi), \\ \notag
|T^2_{reg}| & \geq & 2 |C^2_{reg}|.
\end{eqnarray}

The first inequality is because 
a regular leaf block $L$ of $F$ must contain a vertex of 
$T^1$ which is not the vertex at which the leaf block is attached to the
rest of $F$ (the cut point). Otherwise
the cycle $s$
can only enter and leave $L$ at the cut point, and so $L$ will be 
a leaf block of $G_0  \cup \supp s = G$, contrary to $G$ being two-connected.

The second inequality is because $\Gamma_P$ is indecomposable at every
point $P$. So at all points $P$ at which $\Phi$ is decomposable
we must have $\dim_k \Gamma_P > \dim_k \Gamma_{0P} = \dim_k \Phi_P$.

The third inequality is because the cycle $s$ must enter 
a two-connected component $F$ in at
least two points in order to avoid $F$ becoming a leaf block in 
$G$. 

Adding the equations (\ref{CliffLigInduk}) obtained by
induction we get:
\[ \sum_{(F,\Phi) \in C} d(\Phi) + \sum_{(F,\Phi) \in C^1} l_{reg}(\Phi) 
+ \sum_{(F,\Phi) \in C} q(\Phi) + 2 |C^2_{reg}|
\geq 2 \sum_{(F,\Phi) \in C} \dim_k \Phi. \]
Now 
\[ d(\Gamma) = \sum_{(F,\Phi) \in C} d(\Phi) + |T|. \]
Using this and the inequalities (\ref{CliffLigUlik}) we obtain
\begin{eqnarray*}
& & d(\Gamma) - |T| + |T^1_{reg}| + |T^1_{nreg}| + |T^2_{nreg}| + |T^2_{reg}| \\
& \geq & 2 \sum_{(F,\Phi) \in C} \dim_k \Phi \geq 2 \dim \Gamma_0. 
\end{eqnarray*}
Since $T$ is the union of $T^1_{reg}, T^1_{nreg}, T^2_{reg}$ and $T^2_{nreg}$
we obtain 
\[ d(\Gamma) \geq 2 \dim_k \Gamma_0 \]
and since $\dim_k \Gamma_0 = \dim_k \Gamma - 1$ we get




\[ d(\Gamma) + 2 \geq 2 \dim_k \Gamma. \]
\end{proof}



\begin{remark} \label{Graph_curves}
An interesting topic for further work was pointed out by F.O. Schreyer. When the graph is trivalent, the multiplication in the canonical module gives rise to a line configuration in $\mathbb{P}^{g-1}$ in the following way. The projectification of the surjective multiplication maps $\widetilde{H}_{d-1}(G)\cong(\omega_{k[G]})_0\overset{x_i}{\rightarrow}(\omega_{k[G]})_{\mathbf{e}_i}\cong \widetilde{H}_{d-1}(G ;G-\{i\})$ is a line in $\mathbb{P}^{g-1}$. The line configuration of the union of all these lines is a graph curve of arithmetic genus $g$ which is canonically embedded. Furthermore, a $g^1_d$ for the graph gives a projection of this curve to a $\mathbb{P}^1$ of degree $d$. It is therefore natural to ask if the gonality of the graph and this curve is the same. Similar curves can also be made from more general $2$-connected graphs, by choosing appropriate rational curves in the corresponding union of projective spaces. One interesting question is if there is a connection between the gonality (or an appropriately defined Clifford index) of the graph and the minimal free resolution of this canonical curve as in Green's conjecture for a canonical curve, which says that you can read the Clifford index of a curve of the graded betti numbers of the minimal free resolution. In many examples we made, this seems to be the case, but it seems very difficult to show this in general. This is first of all because the ideal of the canonical curve is not multigraded, and there is not a combinatorial way to describe the minimal free resolution. Secondly, the Clifford index of a graph is also not easily computed.

\end{remark}

\bibliography{hlohne}

\end{document}